\documentclass[12pt]{article}

\title{Transcendence degrees of fields generated by exponentials of products}

\author{Heinrich Massold}

\usepackage{bezier,amsfonts}

\usepackage{amsmath,amsfonts,amssymb}

\newtheorem{Satz}{}[section]

\newcommand{\satz}[1]{\vspace{2mm} \begin{Satz}{\bf #1}}

\newcommand{\proof}{\vspace{4mm} {\sc Proof}\hspace{0.2cm}}

\newcommand{\R}{{\mathbb R}}

\newcommand{\CA}{{\cal A}}

\newcommand{\Z}{\mathbb{Z}}

\newcommand{\N}{\mathbb{N}}

\newcommand{\CH}{{\cal H}}

\newcommand{\CD}{{\cal D}}

\newcommand{\CF}{{\cal F}}

\newcommand{\Q}{{\mathbb Q}}

\newcommand{\C}{{\mathbb C}}

\newcommand{\G}{{\mathbb G}}

\newcommand{\z}{{\bf z}}

\newcommand{\x}{{\bf x}}

\newcommand{\gen}{\mbox{gen}}

\begin{document}

\parindent0mm

\maketitle

\thispagestyle{empty}

\begin{abstract}
Let $\theta=(\theta_1,\ldots,\theta_m) \in \R^m, \kappa=(\kappa_1,\ldots,\kappa_n) \in \R^n$ be two
tuples of real numbers each linearly independent over $\Q$, and $T$ the transcendence degree
of the field generated by $\{\exp(\theta_i \kappa_j) | i=1,\ldots,m, \; j=1,\ldots,n \}$ over $\Q$. The
estimate $T \geq \frac{mn}{m+n} -1$ has been conjectured for some time but could only be proved
under additional hypotheses for $\theta$ and $\kappa$. This paper proves a weaker estimate for $T$ while
also reducing the strong estimate to a prominent conjecture on intersections of subvarieties of split tori
with subgroups.
\end{abstract}

\tableofcontents

\section{Introduction}

Let $l$ be a natural number, $G$ the algebraic group $G= (\G_m)_{\Z}^l$ and
$X \subset G$ an irreducible subvariety defined over an algebraic extension of $\Q$
that is contained in no
proper algebraic subgroup of $G$. With $t = \dim X$, an irreducible 
subvariety $Y \subset X$ of dimension $s$ is called special if $Y$ is
contained in an algebraic subgroup $H \subset G$ of codimension $t+1-s$.

\satz{Conjecture} \label{ceins}
In the situation above, there is a proper Zariski closed subset
$\bar{X} \subset X$ that contains all special subvarieties.
\end{Satz}

This conjecture is proved for $t=1$ (\cite{Mau}, Th\'eor\`eme 1.2).

\vspace{2mm}


For a finite set
$\Theta = \{\Theta_1,\ldots,\Theta_n\} \subset \C$, denote the transcendence
degree of $\Q(\Theta)$ over $\Q$ by $T(\Theta)$.

A tuple of complex numbers $(\theta_1,\ldots,\theta_m) \in \C^m$ will be called regular if the point
$(\exp (\theta_1),\ldots,\exp(\theta_m)) \in \G_m^m(\C)$ is not contained in any proper algebraic 
subgroup, which is equivalent to
$\theta_1,\ldots,\theta_m, \pi i$ being linearly independent over $\Q$.

If $\theta_1,\ldots, \theta_m  \in \R$  are linearly independent over $\Q$, then 
$(\theta_1,\ldots,\theta_m)$ is regular. If 
$\theta_1, \ldots, \theta_m  \in \C$ are linearly independent over $\Q$ but $(\theta_1,\ldots,\theta_m)$ is
not regular, then there is a subset
$\{ \theta_{i_1},\ldots,\theta_{i_{m-1}} \}\subset \{ \theta_1,\ldots,\theta_m \in \C \}$ of $m-1$ 
numbers such that $\theta_{i_1},\ldots,\theta_{i_{m-1}}, \pi i$ are linearly independent over $\Q$ and
\[ T(\exp(\theta_{i_1}),\ldots,\exp(\theta_{i_{m-1}})) = T(\exp(\theta_1),\ldots,\exp(\theta_m)). \]

\satz{Definition}
Let $n,m,\nu,\mu \in \N, \nu \leq n, \mu \leq m, \eta>1$, and
$\theta=(\theta_1, \ldots,\theta_m)\in \C^m$ and $\kappa=(\kappa_1,\ldots,\kappa_n) \in \C^n$
be regular.

\begin{enumerate}

\item

The $m$-tuple
$\theta=\{ \theta_1, \ldots,\theta_m \}$ is called $(\mu,\eta)$-generic, if there is a nonsingular 
$m \times m$-matrix $A$ with entries in $\Q$ such that for every sufficiently big 
$D \in \N$, there are $1 \leq i_1< \cdots <i_\mu \leq m$ such that for every nonzero
$l=(l_1,\ldots,l_\mu) \in \Z^\mu$ with $|l| \leq D$, the inequality
\[ \log |\exp((l_1A \theta)_{i_1} + \cdots + l_\mu (A \theta)_{i_\mu})-1| \gg  - D^\eta \]
holds, where the implied constant depends only on $\theta$ and $A$.

The $m$-tuple $\theta$ is called $(\mu,\eta)$-special if it is not 
$(\mu,\eta)$-generic, i.e.\@ if for every $c>0$ and every nonsingular
$A \in M_{m \times m}(\Q)$, there are infinitely many $D \in \N$ 
such that for every $\{i_1,\ldots,i_\mu\} \subset \{1,\ldots, m\}$ there is a nonzero
$l=(l_1, \ldots, l_\mu) \in \Z^\mu$ with $|l|\leq D$ and
\[ |\exp(l_1 (A \theta)_{i_1} + \cdots l_\mu (A \theta)_{i_\mu})-1| \leq - cD^\eta. \]

Denote by $\mbox{gen}(\theta,\eta)$ 
the biggest number $\nu \in \N$ such that $\theta$ is $(\nu,\eta)$-generic.

\item
The bituple
$((\theta_1, \ldots,\theta_m),(\kappa_1,\ldots,\kappa_n)) $ is called $(\mu,\nu,\eta)$-generic
if there are regular matrices $A \in M_{m \times m}(\Q),B \in M_{n \times n}
(\Q)$
and a $D_0 \in \N$ such that for all $L,R \geq D_0$, there are
subsets $\{i_1,\ldots,i_\mu \} \subset \{1,\ldots,m\}$,
$\{j_1,\ldots,j_\nu \} \subset \{1,\ldots,n\}$ such that for all nonzero
$l=(l_1, \ldots, l_\mu) \in \Z^\mu,  r=(r_1,\ldots,r_\nu) \in \Z^\nu$ with 
$|l| \leq L$ and $|r| \leq r$, the inequality
\[ \log |\exp((l_1 (A \theta)_{i_1} + \cdots l_\mu (A \theta)_{i_\mu})
   (r_1 (B \kappa)_{j_1} + \cdots r_\nu (B \kappa)_{j_\nu}))-1|\gg  - L^\eta - R^\eta \]
holds, where the implied constant depends only on $\theta,\kappa, A$, and $B$. 

The bituple
$((\theta_1, \ldots,\theta_m),(\kappa_1,\ldots,\kappa_n)) $ is called $(\mu,\nu,\eta)$-special if it is not \\
$(\mu,\nu,\eta)$-generic, i.e.\@ if for every $c>0$ and every 
$A \in M_{n \times n}(\Q),B \in M_{n \times n}(\Q)$ with nonzero determinant, there are
arbitrarily big $L,R \in \N$ 
such that for every $\{i_1,\ldots,i_\mu\} \subset \{1,\ldots,n\}$ there are nonzero
$l=(l_1, \ldots, l_\mu) \in \Z^\mu, r=(r_1,\ldots,r_\nu) \in \Z^\mu$ with
 $|l|\leq L, |r| \leq R$, and
\[ \log |\exp((l_1 (A \theta)_{i_1} + \cdots l_\mu (A \theta)_{i_\mu})
   (r_1 (B \kappa)_{j_1} + \cdots r_\nu (B \kappa)_{j_\nu}))-1| \leq  -c( L^\eta +R^\eta) \]


\end{enumerate}

\end{Satz}


\satz{Conjecture} \label{czwei}
For $m,n \in\N$, let $\theta=(\theta_1,\ldots,\theta_n) \in \C^m$ and
$\kappa=(\kappa_1,\ldots,\kappa_n) \in \in \C^n$ be regular $\Q$.

\begin{enumerate}

\item
If $\gen(\theta,t) \leq m-t$ for some $t \in \N$, then $T(\exp(\theta_1),\ldots,\exp(\theta_m)) \geq t$, 
and if
$\gen(\kappa,t) \leq n-t$ for some $t \in \N$, then $T(\exp(\kappa_1),\ldots,\exp(\kappa_n)) \geq t$.

\item
If either $\mbox{gen} (\theta,t) < m$, and $n\geq t$, or
$\mbox{gen}(\kappa,t) <n$, and $m \geq t$ for some $t \in \N$, then
\[ T(\{ \exp(\theta_i \kappa_j)| i=1,\ldots,m, \; j=1,\ldots,n\}) \geq t. \]

\item
If $\theta_1,\ldots,\theta_n, \kappa_1,\ldots,\kappa_m$ are all real, then
\[ T(\{ \exp(\theta_i \kappa_j)| i=1,\ldots,m, \; j=1,\ldots,n\}) \geq\frac{mn}{m+n}-1, \]
for every $\epsilon >0$.
\end{enumerate}

\end{Satz}

\satz{Theorem I}
Conjecture \ref{ceins} implies conjecture \ref{czwei}
\end{Satz}


\satz{Theorem II}
For $m,n \in \N$, let $\theta=(\theta_1,\ldots,\theta_m) \in \C1m$ and
$\kappa_1,\ldots,\kappa_n \in \C^n$ be regular.
\begin{enumerate}

\item
If $\mbox{gen} (\theta,t) \leq \frac{m-t}{\max(t-1,1)}$ for some $t \in \N$, 
then $T(\exp(\theta_1),\ldots,\exp(\theta_m)) \geq t$.
Likewise, if $\mbox{gen} (\kappa,t) \leq \frac{n-t}{\max(t-1,t)}$ for some $t \in \N$, 
then $T(\exp(\kappa_1),\ldots,\exp(\kappa_n)) \geq t$.

\item
For $\theta_1,\ldots,\theta_n, \kappa_1,\ldots,\kappa_m$ all real, $t,m,n  \in \N$, if
there exist $\mu, \nu \in \N$ with $\mu \leq m$, and $\nu \leq n$, 
such that
\[ \frac{\mu \nu}{\mu + \nu} > t, \quad
   \mu \leq \frac {m-t}{\max(t-1,1)}, \quad \mbox{and} \quad
   \nu \leq \frac {n-t}{\max(t-1,1)}, \]
then
\[ T(\{ \exp(\theta_i \kappa_j)| i=1,\ldots,m, \; j=1,\ldots,n\}) \geq t. \]
The above three inequalities are fulfilled e.\@ g.\@ for $t\geq 2$ arbitrary,
$\mu=2t+1, \nu=2t, m \geq 2t^2-1$, and $n \geq 2 t^2-t$.

\end{enumerate}

\end{Satz}

\satz{Corollary}

\item
If $\theta_1,\ldots,\theta_m \in \R$ as well as $\kappa_1,\ldots,\kappa_n \in \R$ are linearly 
independent over $\Q$, then
\[ T(\{ \exp(\theta_i \kappa_j)| i=1,\ldots,m, \; j=1,\ldots,n\}) \geq 
\left[\sqrt{\frac {min(m,n)+1}2}\right]. \]

\item 
If $\zeta \in \R$ is a transcendental number, and $m,n \in \Z$, then
\[ T(\{\exp(\zeta^m),\exp(\zeta^{m+1}),\ldots,\exp(\zeta^{n-m})\}) \geq
   \left[\sqrt{\frac{n}4 +\frac 12} \right]. \]
\end{Satz}

\proof
1. Assume $\min(m,n)=n$. In Theorem II.2 take $t=\left[ \sqrt{\frac {n+1}2}\right]$,
$\mu=\nu= 2t^2 -1$. An easy calculation shows that 
$m \geq n \geq (\nu-1)(t-1) -1$, and $\frac{\mu \nu}{\mu + \nu} > t$, and the estimate
$T(\{ \exp(\theta_i \kappa_j)| i=1,\ldots,m, \; j=1,\ldots,n\}) \geq t =
\left[\sqrt{\frac {n+1} 2}\right]$ follows.

\vspace{1mm}

2. In part one, take $\theta= (\zeta^{\left[ \frac m2\right]}, \ldots,\zeta^{(n-m)-\left[ \frac{n-m}2 \right]})$, and
$\kappa= (\zeta^{m-[\frac m2]},\ldots,\zeta^{\left[ \frac{n-m}2 \right]})$.

\vspace{2mm}

Throughout the whole paper, the norm of a polynomial with coefficients in $\Z$ will always be the
maximum norm, e.\@ g.\@  for $f= \sum_{i=0}^n a_i x^i$, we have $|f| =  \max_{i=0,\ldots,n} |a_i|$.

To prove estimates for transcendence degrees, two special cases of the Philippon criterion will be needed.

\satz{Proposition} \label{phil}
\begin{enumerate}
Let $\Theta=(\Theta_1, \ldots,\Theta_n) \in \C^n$, and $c_1,c_2,\eta>0,c_3 \in \R$,
and denote by $|z_1,z_2|$ the distance of the point $z_1,z_2$ in $\C^n$. 
There is a constant $C \gg 0$, only depending on $n,\Theta,c_1,c_2$ and $\eta$ such that 

\item
if for every sufficiently big natural number $D$, 
there are  polynomials $f_1,\ldots,f_m$  in $n$ variables, such that
\[ \deg f_i \leq c_1 D, \quad \log |f_i| \leq c_2 D, \quad
   \log |f_i(\Theta_1,\ldots,\Theta_n)| \leq -C D^\eta, \quad i=1,\ldots,m, \]
and for every common zero $z=(z_1,\ldots,z_n)$ of $f_1,\ldots,f_m$, we have
$\log |z,\Theta| \geq -3C D^\eta$. Then, $T(\Theta) > \eta-1$.

\item 
if for infinitely many natural numbers $D$,
there are  polynomials $f_1,\ldots,f_m$  in $n$ variables, such that
\[ \deg f_i \leq c_1 D, \quad \log |f_i| \leq c_2 D, \quad
   \log |f_i(\Theta_1,\ldots,\Theta_n)| \leq -C D^\eta, \quad i=1,\ldots,m, \]
and for every common zero $z=(z_1,\ldots,z_n)$ of $f_1,\ldots,f_m$, we have
$ \log |z,\Theta| \geq c_3$. Then, $T(\Theta) > \eta-1$.
\end{enumerate}
\end{Satz}

\proof This follows from \cite{Ph}, Th\'eor\`eme 2.11.

\section{Generic points}

\satz{Lemma} \label{triv}
\begin{enumerate}

\item

For all $\theta,\eta$, we have $1 \leq \mbox{gen}(\theta,\eta) \leq n$.

\item
For a nonsingular $A \in M_{m \times m}(\Q)$,
the tuple $(\theta_1,\ldots,\theta_m)$ is $(\mu,\eta)$-regular if and only
$(A \theta_1,\ldots,A \theta_m)$ is.

\item
For a nonzero $a \in \R$
the tuple $(\theta_1,\ldots,\theta_m) \in \R^m$ is $(\mu,\eta)$-regular if and only
$(a\theta_1,\ldots,a \theta_m)$ is.

\item
For nonsingular $A \in M_{m \times m}(\Q),B\in M_{n \times n}(\Q)$,
the bituple $(\theta,\kappa)$ is $(\mu,\nu,\eta)$-regular if and only
$(A \theta, B\kappa)$ is.

\item
For $A \in M_{s \times m}(\Q)$ a matrix of rank $s$,
\[ \gen(A \theta,t-1) \leq \gen(\theta,t-1) \leq \gen(\theta,t). \]


\item
For $\theta_1,\ldots,\theta_m,\kappa_1, \ldots, \kappa_n \in \R$, if
$\theta=(\theta_1,\ldots,\theta_m)$ is $(\mu,\eta)$-regular and 
$\kappa=(\kappa_1,\ldots, \kappa_n)$ is $(\nu,\eta)$-regular, then the bituple
$(\theta,\kappa)$ is $(\mu,\nu,\eta)$-regular.


\end{enumerate}
\end{Satz}

\proof
1. Obvious.

\vspace{2mm}

2,.3.,4. These claims hold since the relation $\gg$ remains true if one side is changed by a fixed
multiplicative constant.

\vspace{2mm}

5. Let $\mu = \gen(A \theta,t-1)$. Then for every sufficiently big $D$ and every
$I =(i_1,\ldots, i_\mu)$ with $1 \leq i_1 < \cdots < i_\mu \leq s$
and every $l=(l_1, \ldots,l_\mu)$ with $|l| \leq D$, the inequality
\[ \log |l_1(A\theta)_{i_1} +  \cdots + l_\mu(A\theta)_{i_\mu}-1| \gg -D^t. \]
holds. If $A$ is extended to an $(m \times m)$ matrix $A'$ of rank $n$, then
\[ \log |l_1 (A'\theta)_{i_1} + \cdots + l_\mu (A'\theta)_{i_\mu}-1| \gg -D^t, \]
for all $I =(i_1,\ldots,i_\mu)$ with  $1 \leq i_1 < \cdots < i_\mu \leq s$ and all
$l$, since the left hand side is unchanged because of $i_\mu \leq s$.
Hence, $\gen(A \theta,t-1)=\mu \leq \gen(A'\theta,t-1)$. Further,
$\gen(A'\theta,t-1) = \gen(\theta,t-1)$ by part 2 of the Lemma, and the inequality
$\gen(\theta,t-1) \leq \gen(\theta,t)$ trivially holds.

\vspace{2mm}

6. Since $\exp: \R \to \R^+$ is a bijection with $\exp(0) = \exp'(0) =1$, the relations
\[ \log |l_1 \exp((A \theta)_{i_1} + \cdots + l_\mu (A \theta)_{i_\mu})-1| \gg  - L^\eta, \]
\[ \log |l_1 \exp((B \kappa)_{j_1} + \cdots + l_\mu (B \kappa)_{j_\nu})-1| \gg  - R^\eta \]
are equivalent to
\[ \log |l_1 (A \theta)_{i_1} + \cdots + l_\mu (A \theta)_{i_\mu}| \gg  - L^\eta, \]
\[ \log |l_1 (B \kappa)_{j_1} + \cdots + l_\mu (B \kappa)_{j_\nu}| \gg  - R^\eta, \]
which in turn implies
\[ \log |((l_1 (A \theta)_{i_1} + \cdots l_\mu (A \theta)_{i_\mu})
   (r_1 (B \kappa)_{j_1} + \cdots r_\nu (B \kappa)_{j_\nu})|\gg  - L^\eta - R^\eta, \]
which again is equivalent to
\[ \log |\exp((l_1 (A \theta)_{i_1} + \cdots l_\mu (A \theta)_{i_\mu})
   (r_1 (B \kappa)_{j_1} + \cdots r_\nu (B \kappa)_{j_\nu}))-1|\gg  - L^\eta - R^\eta. \]

\vspace{2mm}

\satz{Proposition} \label{gen}
If $\theta=(\theta_1,\ldots,\theta_m)$ as well as $\kappa=(\kappa_1,\ldots,\kappa_n)$ 
are tuples of linearly independent over $\Q$,  and for some $\eta > 1$, the inequality
\[ \frac{\mbox{gen}((\theta),\eta) \; \mbox{gen}((\kappa),\eta)}
   {\mbox{gen}((\theta),\eta) + \mbox{gen}((\kappa),\eta)} > \eta \]
holds, then $T(\{ \exp(\theta_i \kappa_j)| i=1,\ldots,n, \; j=1,\ldots,m\}) \geq \eta-1$.
\end{Satz}

For $\mbox{gen}(\theta)=m$, and $\mbox{gen}(\kappa) =n$, this is proved in \cite{LN}, chapter 14.3.
The general case can be proved in exactly the same way; the only thing that has to be modified is that 
for every value of the approximation parameter $D$, one works with only $\mu$, respectively $\nu$ 
components of $\theta$ and $\kappa$, which may be different components for every $D$. 
The necessary auxiliary 
polynom for each $D$ is then a polynomial in only these ``active'' variables. As this requires notational and some technical adjustments throughout the whole proof, for the convenience of the reader, I will give a full generalized proof in the appendix.

\section{Special points}



Let $G = \G_m^l$, and $X \subset \G_m$ an irreducible subvariety defined over a finite extension of $\Q$,
that is contained in no proper subgroup of $G$.
For $M$ the character module of $G$, and $N \subset M$ a submodule
define the subgroup
\[ H_N := \bigcap_{\chi \in N} \mbox{ker} \chi \subset G, \]
and for a subvariety $Y \subset X$ define $H_Y$ as the smallest subgroup of $G$ that contains $Y$.

\satz{Definition}
For $s \geq t= \dim X$,
\begin{enumerate}

\item
the variety $X$ is called
$s$-regular, if for every submodule
$N \subset M$ of rank $s$, we have $\dim ((\chi_1,\ldots,\chi_s)(Y)) = t$, where 
$\chi_1,\ldots,\chi_s$ is any basis of $N$, and $(\chi_1,\ldots,\chi_s)$ is the corresponding
map $\G_m^n \to \G_m^s$. 

\item
an irreducible subvariety $Y \subset X$ of dimension $d$ 
is called $s$-special if $\dim H_Y < l-s+d$.

\end{enumerate}
\end{Satz}

Conjecture \ref{ceins} says that there is
a proper Zariski closed subset $\bar{X} \subset X$ that contains
all $t$-special subvarieties of $X$.


\satz{Theorem} \label{sspez}
If $X \subset \G_m^n$ is $s$-regular for some $s\geq t = \dim X$ , 
then there is a proper Zariski closed subset $\bar{X} \subset X$ that
contains all $s$-special subvarieties of $X$.
\end{Satz}

\proof
This result is a mainly technical generalization of 
\cite{Ha}, Corollary 3 and can be found in \cite{Ma}.

\satz{Proposition} \label{transs}
Let $\Theta \in G(\C)$ be a point that is contained in no proper algebraic subgroup,
and $X$ the algebraic closure of $\{\Theta\}$ over $\Q$. 
Suppose that for $t \in \N$, some $c_1 \geq 1 $, an arbitrary  $c>0$, and a proper Zariski closed subset 
$\bar{X} \subset X$,
there is an infinite set $\CD \subset \N$  such that for every $D \in \CD$
there is a submodule $N_D  \subset M$ such that $X \cap H_{N_D} \subset \bar{X}$, and for some basis
$\chi_1,\ldots,\chi_r$ of $N_D$ the inequalities
\[ \deg  \chi_i \leq c_1 D, \quad \mbox{and} \quad 
   \log |\chi_i(\Theta)-1| \leq -c D^t,  \quad i=1,\ldots,r \]
hold. Then, $T(\Theta) \geq t$.
\end{Satz}

\proof
As $\bar{X}(\C)$ does not contain $\Theta$, the distance
$\bar{X}(\C)$ to $\Theta$ is positive. Denote this distance by $c_3$, and
let $\CF =\{f_1,\ldots,f_l \}$ be a set of generators of
the ideal of $X$ in $G$. By shrinking $\CF$, if necessary we may assume
$\deg f_i \leq D$, and $\log |f_i| \leq D$ for all $i=1,\ldots,l$, and every $D \in \CD$.
For $D \in \CD$  let $\CF_D := \CF \cup \{\chi_1-1, \ldots,  \chi_r-1\}$.
Then, by assumption, $\bar{X}$ contains the set of common zeros of $\CF_D$, 
hence every common zero of $\CF_D$ has distance at least
$c_3$ to $\Theta$. Also, for $f \in \CF_D$,
\[ \deg f \leq c_1 D, \quad \log |f| \leq D, \quad \mbox{and} \quad 
    \log |f(\Theta)| \leq - c D^t, \]
the second inequality because for $f$  a character, $|f| \leq 1$,
the last inequality because $\log |\chi_i(\Theta)|-1 \leq -cD^t$ by assumption for all $i=1,\ldots,r$, and
$|f(\Theta)| = 0$ for $f \in \CF$. By Proposition \ref{phil}.2, $T(\Theta) >t-1$, from which
$T(\Theta) \geq t$ follows, since $T(\Theta)$ is integral.

\section{Proof of the main Theorems}

\satz{Lemma} \label{genug}
\begin{enumerate}

\item
Let $V$ be a $\Q$-vector space with basis $\{v_1,\ldots,v_n\}$,
and $\nu$ a natural number less or equal $n$.
For every subset $I \subset \{1,\ldots,n\}$ with $|I| = n-\nu$
choose a non-zero vector $w_I = \sum_{i \in I} a_i v_i$.
Then, the subpsace $W_I$ generated by these $w_I$ has dimension at
least $\nu+1$. More specifically, there are subsets $I_1,\ldots I_{\nu+1}$ such
that $w_{I_j},j=1,\ldots,\nu+1$ are linearly independent.

\item
Let $(\theta_1,\ldots,\theta_m) \in \R^m, (\kappa_1, \ldots, \kappa_n) \in \R^n$ be two tuples of real numbers, linearly independent over $\Q$, and $\bar{G} \subset G = G_m^{mn}$ the smallest algebraic
subgroup that contains $(\exp(\theta_i \kappa_j)_{i \leq m, j \leq n}$.
Further, with $x_{ij}, 1 \leq i \leq m, 1 \leq j \leq n$ the coordinate functions of $\G_m^{mn}$, and
 $l_1, \ldots,l_n \in \Z$, not all zero, let $\chi_j, j=1,\ldots, n$ be the characters
$\chi_j = \prod_{i=1}^m x_{i j}^{l_i}$, and $H = \bar{G} \cap_{j=1}^n \ker \chi_j$. 
Then, the codimension of $H$ in $\bar{G}$ equals n.

\end{enumerate}
\end{Satz}

\proof
1. The $w_{I_j}$ are construed inductively: Take $w_{I_1}$ as any of the $w_I$.
If the vectors $w_{I_1},\ldots,w_{I_j}$ with $j< \nu+1$ are given
let $W_j$ be the space generated by them. As $v_i,i=1,\ldots,n$ are 
a basis of $V$, their rest classes $\bar{v}_i \in V/W_j$ generate
$V/W_j$. Hence there are $n-j$ natural numbers $l_j \leq n$
such that $\bar{v}_{l_j}, j=1,\ldots,n-j$ form a basis
of $V/W_j$. Since $n-j \geq n-\nu$, there is a subset $I_{j+1}$ of
$\{l_j|j=1,\ldots,n-j\}$ with $|I|=n-\nu$,
such that the restclass $\bar{w}_I$ of $w_I$ in $V/W_j$ is nonzero.
Conseqently $w_{I_1},\ldots,w_{I_j},w_{I_{j+1}}$ are linearly independent.

\vspace{2mm}

2. Let $X^{\bar{G}}$ be the set of characters of $\G_m^{mn}$ that are one on $\bar{G}$
and $X^H$ the character module generated by the $\chi_j,j=1,\ldots,n$. If
the codimension of $H$ in $\bar{G}$ were smaller than $n$, the intersection
of $X^{\bar{G}}$ with $X^H$ would be nonzero. So assume that there are
$r_1,\ldots,r_n \in \Z$, not all zero, such that
\[ \chi :=\prod_{j=1}^n \chi_j^{r_j} = \prod_{i=1}^m \prod_{j=1}^n x_{ij}^{l_i r_j} \in X^{\bar{G}}. \]
Then, the resstriction of $\chi$ to $\bar{G}$ would be 1, hence
$\chi((\exp(\theta_i \kappa_j))_{i\leq m, j \leq n}) = 1$. Since the $\theta_i, \kappa_j$ are all real, this is
equivalent to
\[ \sum_{i=1}^m \sum_{j=1}^n (l_i \theta_i)(r_j \kappa_j) =  0 ,\quad \mbox{hence} \quad
   \left(\sum_{i=1}^m l_i \theta_i\right) \left( \sum_{j=1}^n r_j \kappa_j \right)= 0. \]
which implies that either $\sum_{i=1}^m l_i \theta_i$ or $ \sum_{j=1}^n r_j \kappa_j=0$, which in turn,
because of the linearily independent conditions, implies that either
$l_1= \ldots = l_m =0$ or $r_1 = \cdots = r_n =0$ in contradiction to the assumptions.

\proof {\sc of Theorem I}
Let $\Theta=(\Theta_1,\ldots,\Theta_m)=\exp(\theta_1,\ldots,\theta_m)$, and
$x_1,\ldots,x_m$ the coordinate functions of $\G_m^m$.
For $l=(i_1,\ldots,l_m) \in \Z^m$ denote the corresponding character 
$\prod_{i=1}^m x_i^{l_i}$ of $\G_m^m$ by $\chi_l$; it has degree
$\sum_{i=1}^m |l_i| \leq \sqrt{m} |l|$. 

1. For $t=1$, assume that $T(\Theta)=0$, and
$\gen(\theta,1) \leq m-1$. Since $\theta$ is $(m,1)$-special, for every $c>0$ there are
infinitely man $D \in \N$ and nonzero $l=(l_1,\ldots,l_m)$ with
\[ |l| \leq D, \quad \mbox{hence} \quad\deg \chi_l \leq \sqrt{m} |l|, \]
and
\[  |\chi_l(\Theta)-1|= |\chi_l(\Theta_1,\ldots,\Theta_m)-1| \leq
   \exp(- c D). \]
Since $\theta$ is regular, we also have $|\chi_l(\Theta)-1|>0$, and since
$\Theta_i$ is algebraic for every $i=1,\ldots,m$, if $c$ is  sufficiently big,
this contradicts the Liouville inequality.

For $t \geq 2$, let $X$ be the algebraic closure of $\{(\Theta_1,\ldots,\Theta_m)\}$ over $\Q$,
and assume $T(\Theta) = \dim X \leq t-1$.
Since, $\gen(\theta,t) \leq  m-t$,  for every $c>0$, there are
infinitely many $D \in \N$ such that for every
$1 \leq i_1 < \cdots < i_{m+1-t} \leq m$, there is a nonzero $l =(l_{i_1},\ldots,l_{m+1-t}) \in \Z^\mu$ with
$|l| \leq D$, hence $\deg \chi_l \leq \sqrt{m} D$, and
\[ \log |\chi_l(\Theta_{1_1}, \ldots, \Theta_{i_{m+1-t}})-1| \leq -c D^t. \]
By part one of the previous Lemma, the rank of the module $N$ generated by these 
$(l_1,\ldots,l_{m-t}) \in \Z^\mu$ is
$t$. Hence, by conjecture I, there is a Zariski closed subset $\bar{X} \subset X$ such that
$\bar{X}$ contains $X \cap H_N$. Thus the first part of the claim follows from Proposition \ref{transs}.
The second part of the claim is proved analogously.

\vspace{2mm}

2. If $\gen(\theta,t) < m$, then $\theta$ is $(m,t)$-special. Hence, for every $c>0$,
there is an infinite subset $\CD \subset \N$ such that for
every $D \in \CD$, there is an $l=l_1, \ldots, l_m$ such that $\deg \chi_l \leq \sqrt{m}D$, and 
\[ \log |\chi_l(\exp(\theta_1), \ldots,\exp(\theta_m))-1| \leq -c D^{t}. \]
Consequently,
\[ \log |\chi_l(\exp(\theta_1 \kappa_j), \ldots,\exp(\theta_m \kappa_j))-1| \leq -c'c D^{t}, \]
with $j$ any natural number less or equal $n$, and $c'$ only depending on $\kappa$.

Next, with $\Theta = (\exp(\theta_i \kappa_j)_{i \leq m, j \leq n}$, let
$X \subset \G_m^{mn}$ be the algebraic closure of $\{\Theta\}$, and
$\bar{G} \subset \G_m^{mn}$ the smallest subgroup of $G_m^{mn}$ that contains $\theta$.
Then, $\bar{G}$ is isomorphic to $\G_m^a$ for some $a \in \N$,because the 
$\exp(\theta_i \kappa_j)$ are all positive real numbers.
Assume $\dim X = T(\Theta) < t \leq n$.
Since $\chi_j=\prod_{j=1}^m x_{ij}^{l_i}, j=1,\ldots,n$, part 2 of the previous Lemma implies that
the codimension of $H= \bar{G} \cap  \cap_{j=1}^n \ker \chi_j$ in $\bar{G}$ is $n > \dim X$.
By conjecture I, the common zeroes of $N_D$ in $X$ are contained in a fixed proper Zariski closed 
subset $\bar{X}$ of $X$, and the claim follows from Proposition \ref{transs}.
If $\gen(\theta,t) <n$ and $m\geq t$, the claim follows analogously

\vspace{2mm}

3. Let $\theta_1,\ldots,\theta_n,\kappa_m,\ldots,\kappa_m$ be as in the conjecture, and
$\eta= mn/(m+n)$. If $\theta$ is $(m,\eta)$-generic and $\kappa$ is $(n, \eta)$-generic,
then
\[ \frac{\gen(\theta,\eta)\gen(\kappa,\eta)}{\gen(\theta,\eta)+\gen(\kappa,\eta)} =
   \frac{mn}{m+n}
    = \eta > \eta-\epsilon \]
for every $\epsilon >0$, and
Proposition \ref{gen} implies $T(\Theta) \geq \eta-1-\epsilon= \frac{mn}{m+n}-1-\epsilon$,
for every $\epsilon >0$ which in turn implies $T(\Theta) \geq \frac{mn}{m+n}-1$.

If $\theta$ is $(m,\eta)$-special, it is also $(m,t)$-special with $t=[\eta]$, hence
$\gen(\theta,t) < m$. Since $n > \frac{mn}{m+n} \geq t$, part 2 implies
\[ T(\Theta) \geq t = \left[ \frac{mn}{m+n} \right] \geq \frac{mn}{m+n}-1. \] 
If $\kappa$ is $(n,\eta)$-special, $T(\Theta) \geq \frac{mn}{m+n}-1$ follows in
the same way.



\proof {\sc of Theorem II}
1. For $t=1$, the claim coincides with conjecture 1.3 for $t=1$, and the proof of this claim did not use
conjecture I.

\vspace{2mm}

For $t=2$, let $\gen(\theta,2) \leq \frac{m-2}1 +1 = m-2$. Then, by Lemma \ref{triv}.5,
$\gen(\theta,1) \leq \gen(\theta,2) \leq m-2 < m-1$, and by the above 
$T(\Theta_1,\ldots,\Theta_m) \geq 1$. Assume that $T(\Theta_1,\ldots,\Theta_m) = 1$.
As $\theta$ is $(m-1,2)$-special, for every constant $c>0$, there is an infinite set $\CD \subset \N$, such
that for every $D \in \CD$, there are $l,\bar{l} \in\Z^m$, with $|l|, |\bar{l}| \leq D$, and 
\[ \log |\chi_l(\Theta_1, \ldots, \Theta_{m-1})-1| \leq -c D^2, \quad
    \log |\chi_{\bar{l}}(\Theta_2, \ldots, \Theta_m)-1| \leq -c D^2. \]
Since $l$ and $\bar{l}$ are linearily independent, \cite{Mau}, The\'or\`me 1.2 implies
that the common zeroes of $\chi_l$ and $\chi_{\bar{l}}$ are contained in a fixed
finite subset $\bar{X}$ of $X$, and the first claim follows from proposition \ref{transs}. The second
claim is proved analogously.

\vspace{2mm}

Assume now $t \geq 3$, and the Theorem be true for $t-1$. Let
$\mbox{gen}(\theta,t) \leq \frac{m-t}{t-1}$, and define $s:= m-\mbox{gen}(\theta,t)-1$. Then, 
\[ m-s-1 = \gen(\theta,t) \leq \frac{m-t}{t-1}. \] 
This equation is equivalent to
\[ m-s-1 \leq \frac{s-t+1}{t-2} \quad \Longleftrightarrow \quad
  \gen(\theta,t) \leq \frac{s-(t-1)}{t-2}. \]
Let $\CA$ be the set of $(s \times m)$-matrices of rank $s$ with coefficients in $\Q$.
By Lemma \ref{triv}.5, for every  $A \in \CA$, 
\[  \gen(A\theta,t-1) \leq \mbox{gen}(\theta,t) \leq \frac{s-(t-1)}{t-2}. \]
The induction hypothesis implies $T(\exp( A\theta)) \geq t-1$ for every $A \in \CA$. Since
$T(\Theta) \geq T(\Theta^A)$, we only need to derive a contradiction from
the assumption that $T(\Theta)= t-1$. 
Assume $T(\Theta)=t-1$. Since $T(\Theta^A)= t-1$ for every
$A \in \CA$ the algebraic closure $X$ of $\{\Theta\}$ over $\Q$ is $s$-regular.
Hence, by Theorem \ref{sspez}, there is a proper Zariski closed subset $\bar{X}$ that contains all
$s$-special subvarietes of $X$.

Since $\theta$ is $(m-s,t)$-special, because of $m-s > m-s-1 = \gen(\theta,t)$,
for every $c$ there
are infinitely many $D \in \N$  such that for every $I =\{i_1,\ldots,i_{n-s}\} \subset \{1,\ldots,m\}$,
there is a nonzero $l_I = (l_{i_1},\ldots,l_{i{_n-s}})$ such that
\[ |l_I| \leq D , \quad \mbox{and} \quad
   \log | \chi_{l_I}(\Theta_1, \cdots,\Theta_{n-s})-1| \leq -cD^t.\]
By Lemma \ref{genug}, the $l_I$ generate a submodule of rank $s+1$, and by the above he intersection
of $X$ with $\cap_I \mbox{ker}(\chi_{l_I})$ is contained in $\bar{X}$. Proposition \ref{transs} implies
$T(\exp(\theta)) \geq t$ in contradiction with the assumption.

Of course, the claim about $\kappa$ is proved in the same way.

\vspace{2mm}

2. If $\gen(\theta,t) \geq \mu$ and $\gen(\theta,t) \geq \nu$, then
\[ \frac{\gen(\theta,t)  \gen(\kappa,t)}{\gen(\theta,t)+ \gen(\kappa,t)} \geq
   \frac{\mu \nu}{\mu + \nu} > t, \]
hence for some $\epsilon > 0$,
\[ \frac{\gen(\theta,t)  \gen(\kappa,t)}{\gen(\theta,t)+ \gen(\kappa,t)} > t+\epsilon, \]
and Proposition \ref{gen} implies
$T(\Theta) \geq t+ \epsilon-1$. Since both $t$ and $T(\Theta)$ are natural numbers, this implies
$T(\Theta) \geq t$.

If $\gen(\theta,t) < \mu \leq \frac{m-t}{\max(1,t-1)}$, by Lemma \ref{triv}.3, likewise
$\gen(\theta_1 \kappa_1, \ldots, \theta_m \kappa_1,t) \leq \frac{m-t}{\max(1,t-1)}$. By part 1,
$T((\exp(\theta_1 \kappa_1), \exp(\theta_m  \kappa_1)) \geq t$. The claim thus follows from the trivial fact
$T(\Theta) \geq T((\exp(\theta_1 \kappa_1), \exp(\theta_m  \kappa_1))$. 
If $\gen(\kappa,t) < \nu$, the claim is proved analogously.

\begin{appendix}

\section{Proof of Proposition \ref{gen}}

Let $\Theta =(\exp(\theta_i \kappa_j))_{i\leq m, j \leq n} \in \G_m^{mn}$, and
\[ \bar{\Theta}_k := (\exp(\theta_i \kappa_j))_{i \leq m, j \leq n, a \leq k}  \in \G_m^{m nk}(\C), \]
where $k$ is any natural number. 
Of course, $T(\bar{\Theta}_k) = T(\Theta)$, for any $k$.

For $\eta$ as in the Proposition, let
$\mu := gen(\theta,\eta), \nu:= gen(\kappa,\eta)$, and for $D\in \N$, let
\[ L=L(D):= [D^{\frac \nu{\mu+\nu}}], \quad R=R(D):= [(2 \mu +1) D^{\frac \mu{\mu+\nu}}]. \]

Further, for $I =(i_1,\ldots,i_\mu), J= (j_1,\ldots,m_\nu)$ with
$1 \leq i_1 < \cdots < i_\mu \leq m$, and $1 \leq j_1 < \cdots < i_\nu \leq n$, let 
\[ p_I: \G_m^{mk} \to G_{I,k} \cong \G_m^{\mu k}, 
   \quad (z_{ia})_{i \leq m, a\leq k} \mapsto (z_{i_\lambda a})_{\lambda \leq \mu, a\leq k}, \]
\[ p_J: \G_m^{nk} \to G_{J,k} \cong 
   \G_m^{\nu k}, \quad (z_{jb})_{j \leq n, b \leq k} \mapsto (z_{j_\rho b})_{\rho \leq \nu, b \leq k}, \]



which induce maps of coordinate rings
\[ p_I^*: \Z[(\G_m)_{I,k}] \to \Z[\G_m^{mk}], \quad
   p_J^*: \Z[(\G_m)_{J,k}] \to \Z[\G_m^{nmk}], \]
For $R \in \N$ let 
$B_R:=\{ r \in M_{J \times k}(\N)| r_{j_\rho a} \leq R, \; \forall \rho=1,\ldots,\mu, a=1,\ldots, k\}$
be the matrices indexed by $J =(j_1,\ldots,j_\nu)$ and $1,\ldots,k$ with entries in the natural numbers less 
or equal to $R$. For $r \in B_R$, and $a \in \{1,\ldots,k\}$  define
\[ m_{r,a}: \G_m^{mn k}\to \G_I \cong \G_m^\mu, \quad  
   (z_{i a j})_{i \leq m, a\leq k, j \leq n} 
   \mapsto \left(\prod_{\rho=1}^{\nu} 
   z_{i_\lambda  a j_\rho}^{r_{j_\rho a}}\right)_{\lambda \leq \mu}, \]
\[  m_r: \G_m^{mn k}\to \G_{I,k} \cong  \G_m^{\mu k},  \quad (z_{iaj})_{i \leq m, j \leq n, a \leq k} \mapsto 
    \left(\prod_{\rho=1}^{\nu} 
   z_{i_\lambda  a j_\rho}^{r_{j_\rho a}}\right)_{\lambda \leq \mu, a \leq k}, \]
with corresponding coordinate maps
\[ m_{r,q}^*: \Z[\G_{I}] \to \Z[\G_m^{m n k}], \quad m^*_r: \Z[\G_{I,k}] \to \Z[\G_m^{m n k}]. \]
We have 
\begin{equation} \label{gradab} 
\deg m^*_{r,a}(f) \leq \max_{\rho=1,\ldots,\nu} |r_{j_\rho a}| \deg f, \quad
   \deg m^*_r(g)\leq \max_{\rho=1,\ldots,\nu,a=1,\ldots,k} |r_{j_\rho a}| \deg g, 
\end{equation}
and 
\begin{equation} \label{hoeab}
  \log|m^*_{r,a}(f)| \leq \log|f|, \quad
    \log|m^*_r (g)| \leq \log |g|, 
\end{equation}
for all $f \in \Z[\G_I], g \in \Z[\G_{I,k}]$.

Finally, for $\z =(z_{iaj})_{iaj} \in \G_m^{mnk}$, $a\leq k$, and $R \in \N$, define
\[ \Sigma_{R,a}(\z) := \{ m_{r,a}(\z) | \; r \in B_R\} \subset G_I(\C), \quad a=1,\ldots,k, \]
\[ \Sigma_R(\z) := \{ m_r (\z) | \; r \in B_R  \}=
   \Sigma_{R,1}(\z) \times \cdots \times \Sigma_{R,k}(\z) \subset G_{I,k}(\C). \]

\satz{Propostion} \label{grmult}
For $G = \G_m^m$, and $\bar{\Sigma} \subset G(\C)$ a finite subset, define
\[ \Sigma = \Sigma(d) = \{ \sigma_1 \cdot \cdots \cdot \sigma_d| \sigma_i \in \bar{\sigma}, i = 1, \ldots,d \}, \]
and assume there is a polynomial
$f \in \C[z_1,\ldots,z_m]$ of degree at most $L$ that vanishes on every point of $\Sigma$.
Then there is a proper subgroup $H \subset G$ such that
\[ \mbox{card} (\Sigma H/H) \CH_H(L)  \leq \CH_G(L), \]
where $\CH_H, \CH_G$ are the Hilbert functions of $H$ and $G$. Moreover $H$ can be chosen as the subgroup given by an equation $\prod_{i=1}^m x_i^{l_i} = 1$
where $|l_i| \leq L$ for every $i=1,\ldots,m$.
\end{Satz}

\proof
In \cite{LN}, ch.\@ 11, Theorem 4.1 take $T=0$.

\vspace{2mm}
\satz{Lemma} \label{prod}
For $l \in \N$, and $\Sigma$ a finite subset of $\C^l$ define
\[ \omega (\Sigma) := \min \{ \deg P| P \in \C[z_1,\ldots,z_l], P \neq 0, P(\sigma)=0 \; 
   \forall \sigma \in \Sigma. \}\]
Then, for $\Sigma_1,\ldots,\Sigma_k$ finite subsets of $\C^\mu$,
\[ \omega(\Sigma_1 \times \cdots \times \Sigma_k) = \min_{1 \leq a \leq k} \omega(\Sigma_a). \]
\end{Satz}

\proof \cite{LN}, ch.\@ 14, Proposition 3.3.

\satz{Lemma} \label{abst}
For $\eta>1$ assume $(\theta,\kappa)$ is $(\nu,\mu,\eta)$-regular, and
$\z = (z_{i a j})_{i a j} = (\z_a)_{a \leq k} \in \G^{mn k}(\C)$ is any point. For every sufficiently big
$D$, there are subsets $I=\{i_1,\ldots,i_\nu\} \subset \{1,\ldots,n\}$ and
$J=\{j_1,\ldots,i_\mu\} \subset \{1,\ldots,m\}$ such that the existence
of a polynomial $f \in \Z[(\G)_I]$ with $\deg f \leq L(D)$ that is
zero on every point of $\Sigma_{\left[\frac{R(D)}2\right]}(\z)$, implies $\log |\Theta_k,\z|\gg -  D^\eta$,
where $|\bar{\Theta}_k, \z|$ is the distance of $\bar{\Theta}_k$ to $\z$, and the implied constant depends
only on $m,n,k$ and $\Theta$.

\end{Satz}

\proof
Let $(l_{\lambda})_{\lambda } \in \Z^{\mu}\setminus \{0\}, (r_{\rho})_{\rho} \in \Z^{\nu}\setminus \{0\}$,
$I \subset \{1,\ldots,m\}$, $J \subset \{1,\ldots,n\}$, and $a \leq k$. Then,
\[ \prod_{\lambda=1}^\mu \prod_{\rho=1}^\nu
   (\bar{\Theta}_k)_{i_\lambda  a j_\rho}^{l_{\lambda} r_{\rho}}  = 
   \prod_{\lambda=1}^\mu \prod_{\rho=1}^\nu \Theta_{i_\lambda j_\rho}^{l_\lambda r_\rho} =
   \exp\left(\left(\sum_{\lambda=1}^\mu l_\lambda\theta_{i_\lambda} \right)
   \left( \sum_{\rho=1}^\nu 
   r_\rho  \kappa_{j_\rho}   \right) \right). \]
If $|l_{\lambda}| \leq L, \forall \lambda$, and
$|r_{\rho}| \leq R, \forall \rho$, since $L=L(D),R=R(D) \to \infty $ when $D \to \infty$, the $(\nu,\mu,\eta)$-regularity of $\Theta$ implies that for any sufficiently big $D$ there are $I,J$ such that for all 
$(l_{\lambda})_{\lambda} \in \Z^{\mu} \setminus \{0\}, (r_{\rho})_{\rho} \in \Z^{\nu} \setminus \{0\}$ with
$|l_{\lambda}| \leq L, \forall \lambda$ and $|r_{\rho}| \leq R, \forall \rho$, the
inequality
\[ \log \left|\prod_{\lambda=1}^\mu \prod_{\rho=1}^\nu
   (\bar{\Theta}_k)_{i_\lambda  a j_\rho}^{l_{\lambda} r_{\rho}}-1\right| \geq
  \log \left|\sum_{\lambda=1}^\mu l_\lambda\theta_{i_\lambda} \right|
   \left| \sum_{\rho=1}^\nu 
   r_\rho  \kappa_{j_\rho}   \right|  +\log 2 \gg \] 
\begin{equation}\label{regf1} L^\eta- R^\eta \geq - (2k(m+n)D)^\eta \gg - D^\eta
\end{equation}
holds, again because of $\exp(0) = \exp'(0)=1$. (Strictly speaking the inequality holds only if 
the argument of the exponential is sufficiently close to $0$, but these are the only arguments we are 
concerned with.)

Assume now that there is an $\bar{f} \in \Z[(\G)_{I,k}]$ with 
$\deg\bar{f} \leq L$ that is zero at every point of $\Sigma_{\left[\frac{R}{2\mu}\right]}(\z)$. 
By Lemma \ref{prod}, there is an $a \leq k$ and an $f \in \Z[(\G)_I]$ with $\deg f \leq L$
that is zero at every point at 
$\Sigma_{\left[\frac R{2\mu}\right],a}(\z)$.
By Propostion \ref{grmult} there is a subgroup $H$ of $G_I$ such that
\[ \mbox{Card}((\Sigma_{\left[\frac R{2\mu}\right],a}(\z) H)/H) \CH_H(L) \leq \CH_{G_I}(L), \]
Morover, $H$ can be chosen to be defined by an equation
\[ \prod_{\lambda=1}^\mu x_{i_\lambda a}^{l_{\lambda}} = 1 \quad \mbox{with} \quad
   |l_{\lambda}| \leq L, \; \lambda=1,\ldots,\mu, \exists \lambda : l_\lambda \neq 0. \]
As $G_I \cong \G_m^{\mu}$, we have $\CH_{G_I}(L) =L^{\mu}$, hence
$\mbox{Card}((\Sigma_{\left[\frac R{2\mu}\right],a}(\z) H)/H) \leq L^{\mu}$.
As because of $\frac{\mu}{\mu+\nu} > \eta >1$, for sufficiently big $D$,
\[ \mbox{Card}(\Sigma_{\left[\frac R{2\mu}\right],a}(\z)) =\left[\frac R{2\mu}\right]^{\nu} >
    [D^{\frac \mu{\mu+\nu}}]^{\nu} = L^{\mu},  \]
this implies that there are two different points 
$\sigma, \bar{\sigma} \in \Sigma_{\left[\frac R{2 \mu}\right],a}(\z)$ such
that $\sigma \bar{\sigma}^{-1} \in H$. With
$\sigma= m_{r,a}(\z), \bar{\sigma}=m_{\bar{r},a} (\z)$, we have
$\sigma \bar{\sigma}^{-1} = m_{r,a} (\z)  \cdot m_{-\bar{r},a}(\z) $, hence
\begin{equation} \label{regf2}
\prod_{\lambda=1}^\mu \prod_{\rho=1}^\nu 
z_{i_\lambda a j_\rho}^{l_{\lambda}(r_{j_\rho a}-\bar{r}_{j_\rho a})} = 1.
\end{equation}
Since $\sigma \neq \bar{\sigma}$, there is an $a$ and a $\rho$ such that 
$r_{j_\rho a} \neq \bar{r}_{j_\rho a}$, and since 
$\sigma,\bar{\sigma} \in \Sigma_{\left[\frac R{2\mu}\right],a}(\z)$,
the inequalty $|r_{j_\rho a}-\bar{r}_{j_\rho a}| \leq \frac R{\mu}$ holds, for every $\rho=1,\ldots,\nu$.

Assume $\log |\bar{\Theta}_k,\z| \leq -c D^\eta$ with some arbitrarily big constant $c$. 
Since the imaginary part of
$(\bar{\Theta}_k)_{iaj}$ is zero for every $i,a,j$, this implies that there are 
$\x_{iaj} \in \C,  i=1,\ldots, m, a=1,\ldots,k,j=1,\ldots,n$ such that $\z_{iaj} = \exp(x_{iaj})$ and
\[ \log |\mbox{im} \; x_{iaj}| \leq -c D^\eta + \log 2, \quad \forall i=1,\ldots,m, \; a=1,\ldots,k, \; j=1,\ldots,n, \]
consequently because of $\eta>1$, for $D$ sufficiently big
\[ K \log |\mbox{im} \; x_{iaj}| \leq -c' D^\eta + \log 2, \quad \forall i=1,\ldots,m, \; a=1,\ldots,k, \; j=1,\ldots,n, \]
for every number $K \leq 2 \mu \nu LR$. Thus, because of (\ref{regf2}),
\[ \log \left|\prod_{\lambda=1}^\mu \prod_{\rho=1}^\nu
  (\bar{\Theta}_k)_{i_\lambda a j_\rho}^{l_{\lambda}(r_{j_\rho a}-\bar{r}_{j_\rho a})} -1 \right| = \]
\[ \log \left|\prod_{\lambda=1}^\mu \prod_{\rho=1}^\nu
  (\bar{\Theta}_k)_{i_\lambda a j_\rho}^{l_{\lambda}(r_{j_\rho a}-\bar{r}_{j_\rho a})} -  
   \prod_{\lambda=1}^\mu \prod_{\rho=1}^\nu 
  \exp(x_{i_\lambda a j_\rho})^{l_{\lambda}(r_{j_\rho a}-\bar{r}_{j_\rho a})}\right| \leq \]
\[ \log \left| \sum_{\lambda=1} ^\mu  \sum_{\rho=1}^\nu l_\lambda(r_{j_\rho a}-\bar{r}_{j_\rho a}) 
   (\theta_{i_\lambda} \kappa_{j_\rho} - \mbox{Re} \; \x_{i_\lambda a j_\rho})   +
   \sum_{\lambda=1} ^\mu  \sum_{\rho=1}^\nu l_\lambda(r_{j_\rho a}-\bar{r}_{j_\rho a})
   \mbox{im} \; x_{iaj} \right| \leq \]
\[ \max(\log (2\mu \nu LR) + \log  |\bar{\Theta}_k, \z|,- c' D^\eta + \log 2) + \log 2, \]
which by assumption is less or equal, which for sufficiently big $D$ is less or equal
\[ -c''D^\eta,\]
where $c''$ is a constant depending only on $m,n,k$, and$\Theta$ times the arbitrarily big chosen 
constant $c$,
which because of $\eta>1$ for $D \gg 0$ contradicts (\ref{regf1}).

\vspace{3mm}



To proof the proposition, a series of auxiliary polynomials fulfilling the conditions in Proposition \ref{phil}
will be construed. The main tool for this is following Lemma.

\satz{Lemma}
For $r<0$ and an holomorphic function $\varphi:\C^k \to \C$ let
\[ |\varphi|_r := \sup_{|z_i| \leq r, i=1,\ldots k} \varphi(z). \]
Let further $M \in \N$, and $\Delta, U$ be positive real numbers. If $(8U)^{k+1} \leq M \Delta$ and 
$\Delta \leq U$, then for any holomorphic functions $\varphi_1,\ldots,\varphi_M: \C^k \to \C$ with
\[ \sum_{l=1}^M |\varphi_l|_{er} \leq \exp(U), \]
there are numbers $h_1, \ldots,h_M \in \Z$ with $\log |h_l| \leq \Delta, l=1,\ldots M$ such that
the function $\varphi = \sum_{l=1}^M h_l \varphi_M$ satisfies $\log |\varphi|_r \leq -U$. 
\end{Satz}

\proof
\cite{Wal}

\vspace{2mm}

For $X= (x_{i_\lambda a})_{\lambda \leq \mu, j \leq k}$, 
and $A_L$ the set of multidegrees  $d=(d_{\lambda a})_{\lambda \leq \mu, a\leq k}$
denote $|d| := \sum_{\lambda =1}^\mu \sum_{a=1}^k d_{\lambda a} $,
and let
\[ A_L := \left\{  X^d = \left. \prod_{\lambda \leq \mu, a  \leq k}  x_{i_\lambda a}^{d_{\lambda a}}
                                                           \right| \; |d| \leq L  \right\} \]
be the monomials of degree at most $L$ in $\Z[\G^{mk}]$, that lie in
the image of $p_I^* :\Z[G_{I,k}] \to \Z[\G_m^{mk}]$. Further,
\[ i_ \theta : \C^k \to \G_m^{mk}(\C), \quad
   (z_1,\ldots,z_k) \mapsto (\exp(\theta_i z_j))_{i \leq m, j \leq k}, \]
and for $X^d \in A_L$,
\[ \varphi_d : \C^k \to \C, \quad \z  \mapsto (X^d \circ i_\theta)(\z).  \]

Further, put
\[ M := |A_L| = {L+\mu k-1 \choose \mu k} \geq \frac{L^{\mu k}}{(\mu k)!}, \quad
   \Delta :=  D \leq LR, \quad \mbox{and} \quad
   U = \frac{(M \Delta)^{\frac 1{k+1}}}8. \]

As $(8U)^{k+1} =M \Delta, \; \Delta \leq U$, and for all $d \in A_L$, the inequalities
\[ \sup_{|\z_a| \leq eR \nu k |\kappa|,\; a=1,\ldots,k} 
    \log \left( \sum_{d \in A_L} \left| \varphi_d((z_a)_{a\leq k}) \right| \right) +\log M \leq  \]
\[ \sup_{|\z_a| \leq eR \nu k |\kappa|,\; a=1,\ldots,k} 
 \sum_{\lambda \leq \mu, a \leq k} d_{\lambda a} \log |\theta_{i_\lambda} z_a|+ \log |A_L|+ \log M \leq \]
\[e L R m n k |\theta| |\kappa| + 2\log M \leq 
   (mnk(2\mu+1) |\theta| |\kappa|+1) D´\leq U \]
hold. The last two inequalities are valid for sufficiently big $D$ because of $LR \leq (2\mu+1) D$, and
$\log M \ll L \leq D$ for big $D$.
The Lemma above thus implies that for every 
$X^d\in A_L$ there is an $h_d \in \Z$ with $\log |h_d| \leq \Delta$  such that
$\varphi = \sum_{d \in A_L} h_d \varphi_d = \sum_{d \in A_L} h_d (X^d \circ i_\theta)$ 
fulfills
\begin{equation} \label{absch} 
\sup_{|(z_a)_{a\leq k}| \leq R \nu k |\kappa|} \log |\varphi| \leq -U =
   -\frac{(M\Delta)^{\frac 1{k+1}}}8 \leq 
   -\frac {1}{8((\mu k)!)^{\frac 1{k+1}}} D^{\frac{\mu \nu}{\mu + \nu} \frac{k}{k+1}}, 
\end{equation}
where $L= [D^{\frac \nu{\mu + \nu}}]$ was used.
Define $f = \sum_{d \in A_L} X^d$ so that $\varphi = f \circ i_\theta$.

Since for $r \in B_R$, with $r^t$ the transpose of $r$, the inequality
\[ |r^t  \kappa| =
   \left| \left( \sum_{\rho=1}^\nu r_{j_\rho a} \kappa_{j_\rho} \right)_{a \leq k} \right|  \leq
    k \nu R |\kappa| \]
holds, by (\ref{absch})
\[ \log |(m_r^* f) (\bar{\Theta}_k)| = \log |f(m_r (\bar{\Theta}_k))| =
   \log|f((\bar{\Theta}_k)_{i_\lambda a j_\rho}^{r_{j_\rho a}})_{\lambda \leq \mu, a \leq k})| = \]
\[ \log|f((\theta_{i_\lambda} \kappa_{j_\rho})^{r_{j_\rho a}})_{\lambda \leq \mu, a \leq k})|  = 
   \log \left|\varphi\left( \left(\sum_{\rho=1}^\nu   r_{i_\rho a} \kappa_{j_\rho}\right)_{a \leq k} \right)\right| = \]
\[ \log| \varphi( r^t \kappa)|  \leq -c_2 \left(D^{\frac{\mu\nu}{\mu + \nu}}\right)^{\frac k {k+1}}, \]
for all $r \in B_R$.

By the assumption on $\eta$ 
\[ \frac{\mu \nu}{\mu + \nu} = 
   \frac{\mbox{gen}(\theta,\eta) \mbox{gen}(\kappa,\eta)}{\mbox{gen}(\theta,\eta) + \mbox{gen}(\kappa,\eta)}
   > \eta. \]
Thus, for a sufficiently big $k$, also $\left(\frac{\mu \nu}{\mu + \nu}\right)^{\frac k {k+1}} > \eta$ 

Hence for an arbitrarily small $\epsilon >0$ and a sufficiently big $D$,
\begin{equation} \label{noben}
\log|(m_r^*f) (\bar{\Theta}_k)| \leq 
-c_2 \left(D^{\frac{\mu \nu}{\mu + \nu}}\right)^{\frac k{k+1}} \leq
-c_2 D^{\eta} \leq -c D^{\eta-\epsilon}. 
\end{equation}
for all $R \in B_R$, and every $c> 0$.

Assume that there is a point $z  \in \G_m^{mn k}(\C)$ such that
\[ f(m_r (z)) =    ((m_r^* f) (z) = 0, \]
 for all $r \in B_R$. Then, by Lemma \ref{triv}.6 $\Theta$ is $(\mu,\nu,\eta)$-regular, and
 Lemma  \ref{abst} implies
\[ \log |\bar{\Theta}_k, z| \gg - 3D^\eta, \]
hence for a sufficiently big $D$,
\[ \log |\bar{\Theta}_k,z| \geq - c D^{\eta-\epsilon}. \]
Thus, the distance of $\bar{\Theta}_k$ to any common zero of the
set $\{ m_r f | r \in B_R \}$
 is at least $-c D^{\eta-\epsilon}$. Because of (\ref{noben}), to apply Proposition \ref{phil}, it only remains to check the upper bounds on the degree
and length of the $m_r^* f$. By construction $\deg f \leq L, \log |f| \leq LR$. 
Hence, by (\ref{gradab}),
\[ \deg m_r^* f \leq L \max_{\rho_1,\ldots,\nu,a=1,\ldots,a} r_{j_\rho a} \leq LR \leq (2 \mu +1) D, \]
and by (\ref{hoeab})
\[ \log |m_r^* f| \leq LR \leq c'(2\mu+1) D. \]
As all conditions of Proposition \ref{phil} are fulfilled, $T(\bar{\Theta}_k) \geq \eta-\epsilon-1$ follows.
Since $T(\bar{\Theta}_k)$ is integral, 
for a sufficiently small $\epsilon$, this implies $T(\bar{\Theta}_k) \geq \eta-1$, hence
\[ T(\Theta) = T(\bar{\Theta}_k) \geq \eta-1. \]

\end{appendix}


\begin{thebibliography}{0mm}

\bibitem[Ha]{Ha} P.\@ Habegger: On the bounded height conjecture. International Mathematics Research Notices, Volume 2009, Issue 5, 2009, Pages 860–886, 


\bibitem[LNM 1752]{LN} Introduction to algebraic independence theory. Springer 
Lecture Notes 1752


\bibitem[Ma]{Ma} H.\@ Massold: Intersections of subvarietes of $\G_m^n$ with algebraic subgroups. 
To appear

\bibitem[Mau]{Mau} G.\@ Maurin: Courbes alg\'ebrique et \'equations multiplicatives. Mathematische
Annalen, Volume 341, pages 789–824, (2008)

\bibitem[Ph]{Ph} P.\@ Phillippon: Crit\`eres pour l'ind\'ependance alg\'ebriques, Inst. hautes \'etudes sci.\@
publ.\@ math\'ematiques 64 (1986), 5-52

\bibitem[Wal]{Wal} M.\@ Waldschmidt: Transcendance et exponentielles en plusieurs variables, Invent.\@ 
Math.\@ 63/1, (1981), 97-127

\end{thebibliography}
\end{document}